# Modeling the maximum edge-weight *k*-plex partitioning problem

Pedro Martins[1]

Polytechnic Institute of Coimbra – ISCAC, Portugal

Center for Mathematics, Fundamental Applications and Operations Research (CMAF-CIO), Faculty of Sciences, University of Lisbon, Portugal

**Abstract**

Given a sparse undirected graph $G$ with weights on the edges, a $k$-plex partition of $G$ is a partition of its set of nodes such that each component is a $k$-plex. A subset of nodes $S$ is a $k$-plex if the degree of every node in the associated induced subgraph is at least $|S| - k$. The maximum edge-weight $k$-plex partitioning (Max-E$k$PP) problem is to find a $k$-plex partition with maximum total weight, where the partition's weight is the sum of the weights on the edges in the solution.

When $k = 1$, all components in the partition are cliques and the problem becomes the well-known maximum edge-weight clique partitioning (Max-ECP). However, and to our best knowledge, when $k > 1$, the problem has never been modeled. Actually, the literature on the $k$-plex addresses the search for a single component in an unweighted graph.

We propose a polynomial size integer linear programming formulation for the Max-E$k$PP problem and consider the inclusion of additional topological constraints in the model. These constraints involve lower and upper limit capacity bounds in each component and upper bound constraints on the number of components in the final solution. All these characterizations preserve linearity and the initial polynomial size of the model.

We also present computational tests in order to show the models' performance under different parameters' settings. These tests resort to benchmark and real-world graphs.

---

[1] Corresponding author. Address: ISCAC – Quinta Agrícola – Bencanta, 3040-316 Coimbra, Portugal. Tel.: +351 239 802 000; fax: +351 239 445 445.
E-mail address: pmartins@iscac.pt

# 1. Introduction

In this paper we discuss a class of partitioning problems in undirected edge-weighted graphs, where each component in the partition characterizes a $k$-plex. In addition, we also discuss the inclusion of additional topological constraints, involving the number of nodes in each component and the number of components in the final solution. The study is mainly focused on modeling aspects, addressing sparse graphs. We will also devote attention to the applied perspective of the problems.

Let $G = (V,E)$ be a sparse undirected graph, where $V = \{1,\ldots,n\}$ is the set of nodes and $E \subset V^2$ the set of edges, with edge weights $w_{ij} \in \mathbb{R}$, for all $(i,j) \in E$.

For a given integer $k \geq 0$, a $k$-plex in $G$ is a subset of nodes $S \subseteq V$ where the degree of every node in the associated induced subgraph is at least $|S| - k$. When $k = 1$ a $k$-plex represents a clique. However, when $k > 1$, the $k$-plex becomes a degree based relaxation of the clique. This problem has been first proposed in Seidman and Foster (1978), involving the identification of a maximum cardinality $k$-plex in an unweighted sparse graph, known as the maximum $k$-plex problem (Max-$k$P problem). The problem was shown to be NP-hard in Balasundaram et al. (2011). Formulations for the Max-$k$P have been proposed in Balasundaram et al. (2011) and Martins (2010), while heuristics were discussed in Moser et al. (2009) and McClosky and Hicks (2012). Properties of $k$-plexes and a comparison with other clique's relaxation concepts can be found in Pattilo et al. (2013). The Max-$k$P problem has applications, for instance, in HIV transmission (Rothenberg et al. 1998; Rothenberg et al. 2000), in social networks (Mukherjee and Holder 2004), in protein-protein interaction networks (Huber et al. 2007; Martins 2010), in text mining (Balasundaram 2008) and in stock markets (Boginsky et al. 2014). A version that involves weights on the nodes is discussed in Boginsky et al. (2014).

So far, the $k$-plex problem has been addressed ignoring the existence of weights on the edges of the graph. So, we start defining the edge-weight of a $k$-plex as the sum of the weights of all the edges in the induced subgraph. This is a straightforward adaptation of the edge-weight of a clique. In effect, cliques with maximum edge-weight have long been discussed in the literature (see, e.g., Park et al. 1996; Macambira and de Souza 2000; Gouveia and Martins 2015), but the same does not hold for edge-weight $k$-plexes, to our best knowledge. Actually, we can find in the literature a relevant number of applications involving the maximum edge-weight clique problem, namely in protein threading and alignment (Akutsu et al. 2003; Akutsu et al. 2004), in protein side chain packing (Brown et al. 2006) and in market basket analysis (Cavique 2007). However, in many social and biological networks, a clique component can be much restrictive, namely when dealing with missing links or false negative connections in the graph. In those cases, cliques have been criticized for their overly restricted nature, which has motivated the emergence of the *k*-plex concept, among other relaxed versions of cliques (see, e.g., Pattilo et al. 2013).

Unlike the former approaches that seek for a single component, the present paper proposes searching for a partition of *G* into distinct components. In this context, there is a significant research work in the literature addressing the maximum edge-weight cliques partitioning (Max-ECP) problem (see, e.g., Grötschel and Wakabayashi 1989, 1990; Dorndorf and Pesch 1994; Ferreira et al. 1996; Hansen and Jaumard 1997; Mehrotra and Trick 1998, Wang et al. 2006; Oosten et al. 2007; Punnen and Zhang 2012; Sukegawa and Miyauchi 2013; Brimberg et al. 2015; Zhou et al. 2015). However, we have found no references involving an edge-weight *k*-plex partition version, which motivates the effort on the present work. Like the Max-ECP problem, a *k*-plex version can also be seen as a non-hierarchical clustering methodology, contributing to the exploratory class of techniques within data mining. The Max-ECP has been applied, for instance, in clustering and classification problems (Grötschel and Wakabayashi 1989), in stock market analysis (Boginski et al. 2006), in gene expression networks (Kochenberger et al. 2005; Pirim et al. 2014) and in other biological networks (Hüffner et al. 2014).

A closely related problem, involving the partition of a graph into co-*k*-plexes is addressed in Cowen et al. (1997) and Trukhanov et al. (2013). This problem is known in the literature as defective coloring and it also belongs to the NP-hard class. A co-*k*-plex is a subset $S \subseteq V$ where the degree of every node in the associated induced subgraph is at most $k-1$. It can be seen as a complementary concept of a *k*-plex. When $k = 1$, the co-*k*-plex is an independent set.

Most of the mentioned approaches on the Max-ECP assume that the graph is complete, namely when modeling the problem. Following a different direction, we discuss the problem on sparse graphs. We further discuss the incorporation of additional constraints addressing upper and lower limits on the sum of the nodes' weights in each component (when *G* includes weights on the nodes). These constraints have also been discussed in the Max-ECP literature, namely in Johnson et al. (1993); Ferreira et al. (1996); Mehrotra and Trick (1998); Ji and Mitchell (2007) and Oosten et al. (2007).

In the present paper we propose discussing a partitioning problem where each component is a $k$-plex. A component, in this case, can possibly represent an unconnected subgraph. So, given an integer $k \geq 0$, we want to find a $k$-plex partition of $G$ with maximum total weight, denoted as the Max-E$k$PP problem. Considering the relationship among $k$-plexes and cliques, the Max-E$k$PP can be seen as a degree relaxation version of the Max-ECP. In fact, when $k = 1$ the two problems coincide. To our best knowledge, the Max-E$k$PP has never been discussed before. In addition, our motivation stems from the extensive number of applications on both Max-$k$P and the Max-ECP problems and the mentioned overly restrictive nature being endorsed to cliques. We also discuss the inclusion of additional topological constraints, involving upper and lower limits on the sum of the nodes' weights in each component, and an upper limit on the number of components in the final solution.

In this paper, we restrict the discussion to linear and polynomial sized (number of variables and constraints) formulations.

In the next section we provide a detailed description of the Max-E$k$PP problem. In Section 3 we formulate the Max-E$k$PP. Computational tests are conducted in Section 4 and the paper ends with a conclusions section.

## 2. The maximum *k*-plex partitioning problem

Given a simple undirected graph $G = (V,E)$, with $V = \{1,\ldots,n\}$ the set of nodes and $E \subset V^2$ the set of edges, a partition in $G$ is characterized by a partition of its set of nodes $V$. When the partition involves $p$ components, where each component is represented by $V_i$ ($i=1,\ldots,p$), we have $V = V_1 \cup \ldots \cup V_p$ with $V_i \cap V_j = \phi$ for all $i,j=1,\ldots,p$ and $i < j$. The subgraph of $G$ induced by $V_i$ is denoted by $G_i$, for all $i=1,\ldots,p$, characterizing the associated partition in $G$. Each edge $(i,j) \in E$ has an associated weight $w_{ij} \in \mathbb{R}$. The weight of a clique or the weight of a $k$-plex is the sum of all their edge weights. If the graph further includes weights on the nodes, we also define non-negative parameters $q_i$, representing the weight of node $i$, for all $i \in V$.

The neighborhood of node $i$ represents the set of edges incident to $i$ in $G$, being defined by $\delta(i) = \{(i,j) \in E: j \in V\}$, for all $i \in V$, with $d_i = |\delta(i)|$ representing the degree of node $i$ in $G$.

Still within notations, the complementary graph of $G$ is represented by $G^c = (V, E^c)$, with $E^c$ the set of missing edges in $G$, that is, $E^c = \{(i,j) \in V^2 : i > j \text{ and } (i,j) \notin E\}\}$. We also denote by $\delta^c(i)$ the set of edges incident do $i$ in $G^c$.

For a given integer $k \geq 0$, the maximum $k$-plex partitioning (Max-E$k$PP) problem is to find a partition of $G$ with maximum total weight, where each component is a $k$-plex. When $k = 1$, the Max-E$k$PP is the well-known maximum (edge-weight) clique partitioning problem (Max-ECP).

As mentioned in Section 1, the Max-ECP has been discussed assuming that $G$ is a complete graph. However, there is no reason to prevent its usage when the graph is sparse. In effect, a large number of real-world problems are characterized in sparse graphs, where the missing edges truly mean that the connection does not exist. So, it makes sense to discuss the Max-ECP in sparse graphs. On the other hand, when considering the Max-E$k$PP problem with $k > 1$, we should only expect dealing with sparse graphs.

For exemplifying, we consider the 8 nodes sparse graph introduced in Figure 1. The graph has 17 edges and density 0.607.

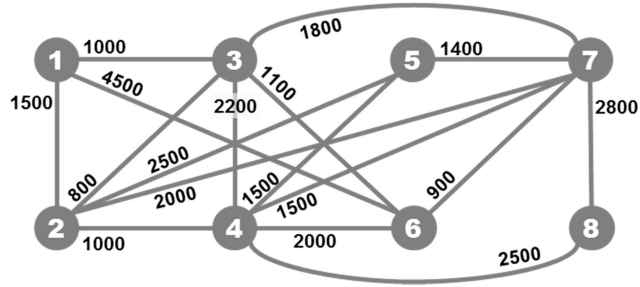

**Figure 1:** Small weighted undirected graph with 8 nodes.

Figure 2 shows the optimum solutions for the Max-E$k$PP problem for $k = 1$, 2 and 3, in images (a), (b) and (c), respectively. When $k = 4$ or 5, the solution involves two components: $V_1 = \{1, 2, 3, 4, 5, 6, 7\}$ and $V_2 = \{8\}$; while for $k \geq 6$ the solution is the entire graph $G$. The optimum values are 15500, 18300, 20800, 24700 and 30000, for $k = 1, 2, 3, 4$ (or 5) and 6, respectively.

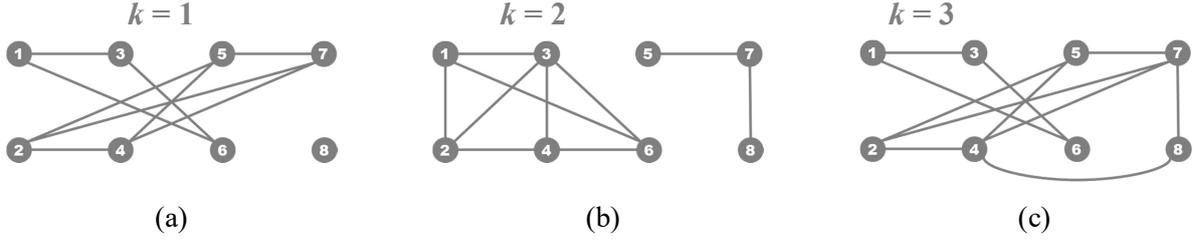

**Figure 2:** Optimum solutions of the Max-E*k*PP problem, for *k* = 1, 2 and 3.

As expected, each node *i* in a component *S* from a Max-E*k*PP solution has degree at least |*S*| − *k*. For instance, in the solution for *k* = 3, node 8 must have degree at least |*S*| − *k* = 5 − 3 = 2 in the subgraph induced by *S*, which means that it is allowed to miss at most *k* − 1 = 2 nodes in component *S*. In fact, for any component *S* in a feasible solution of the Max-E*k*PP, $d_i \geq |S| - k$, for all $i \in S$. So, an immediate result of the *k*-plex definition establishes that if $d_i \geq |V| - k$, for all $i \in V$, then an optimum solution of the Max-E*k*PP is the entire set *V*, assuming that $w_{ij} \geq 0$. This suggests a direct relationship *k* and the minimum degree among all nodes in *G*, when discussing the Max-E*k*PP problem with all edges with non-negative weights, established in the following proposition.

**Proposition 1:** Consider that $w_{ij} \geq 0$ for all $(i,j) \in E$. If $k \geq n - \min_{i \in V}\{d_i\}$, then an Max-E*k*PP optimum solution is $V_1 = V$.

Taken again the example in Figure 1, as $\min_{i \in V}\{d_i\} = 2$, then the optimum solution becomes the entire graph for $k \geq 6$.

When the graph includes edges with negative weights, then the solutions may still involve more than a single component, no matter the value of *k*. To exemplify, consider the graph in Figure 1 with a single change in the weight of edge (1,6), becoming $w_{57} = -4500$. Figure 3 shows the optimum solutions of Max-E*k*PP for *k* = 1, 2 and 3 on the modified graph. The optimum values are 11200, 13800 and 17800, for *k* = 1, 2 and 3, respectively.

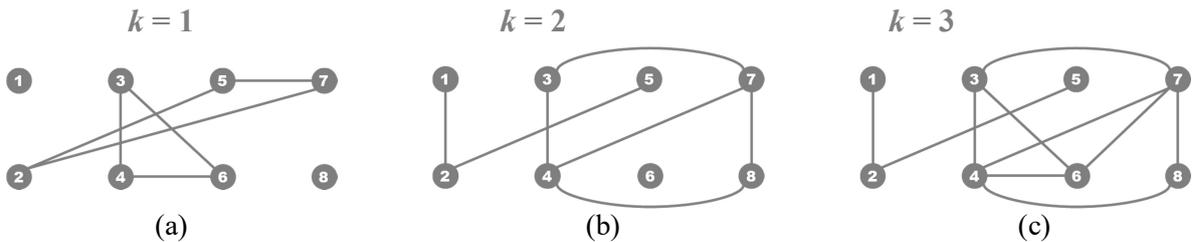

**Figure 3:** Optimum solutions of the Max-E*k*PP problem for *k* = 1, 2, 3, on the graph with $w_{16} = -4500$.

When $k = 4$, the optimum partition is $V_1 = \{2, 3, 4, 5, 7, 8\}$, $V_2 = \{1\}$ and $V_3 = \{6\}$, with total cost 19000. For $k \geq 5$, the optimum partition is $V_1 = \{2, 3, 4, 5, 6, 7, 8\}$ and $V_2 = \{1\}$, with total cost 23000.

In this case, all the solutions avoid the inclusion of edge (1,6) due to its heavy negative weight. Even in the version with $k \geq 6$, in which the degree connectivity condition is entirely relaxed, the optimum solution brought two independent components, showing that Proposition 1 may not hold when the graph includes negative weights on the edges. On the other hand, if the end nodes of an edge with negative weight are in the same component, then the edge must belong to the solution.

Another immediate result of the $k$-plex states that a component $S$ with $|S| \leq k$ may correspond to a set of singletons, contradicting the cohesiveness principle of a component. These components include no edges and we denote them as *spurious*, being treated as isolated nodes. As expected, spurious components are more frequent when $k$ increases and when the graph becomes sparser. Other less cohesive components can also come up when $|S| \leq k$ or $|S|$ is not much larger than $k$. These may include unconnected components. In this study, we do not give much relevancy to those components because they do not represent the main stream of the partition discussion. We will concentrate our attention in the larger and heaviest components instead.

When discussing the Max-E$k$PP problem as a clustering non-hierarchical approach with a given topological structure, we may be interested on the inclusion of additional conditions involving the size and the number of components/clusters in the final solution. Some of those conditions can include upper and lower limits on the sum of the nodes' weights in each component/cluster, and an upper limit on the number of components/clusters in the solution. In order to exemplify, we consider the results returned by the Max-E$k$PP problem for different values of the upper limit on the sum of the nodes' weights in each component (parameter $Q_2$) for $k = 1, 2$ and 3. To this purpose, we take the example in Figure 1 and assume that the nodes' weights are all equal to 1 ($q_i = 1$, for all $i \in V$). Figure 4 shows the results, considering $Q_2 = 2, 3, 4$ and 5.

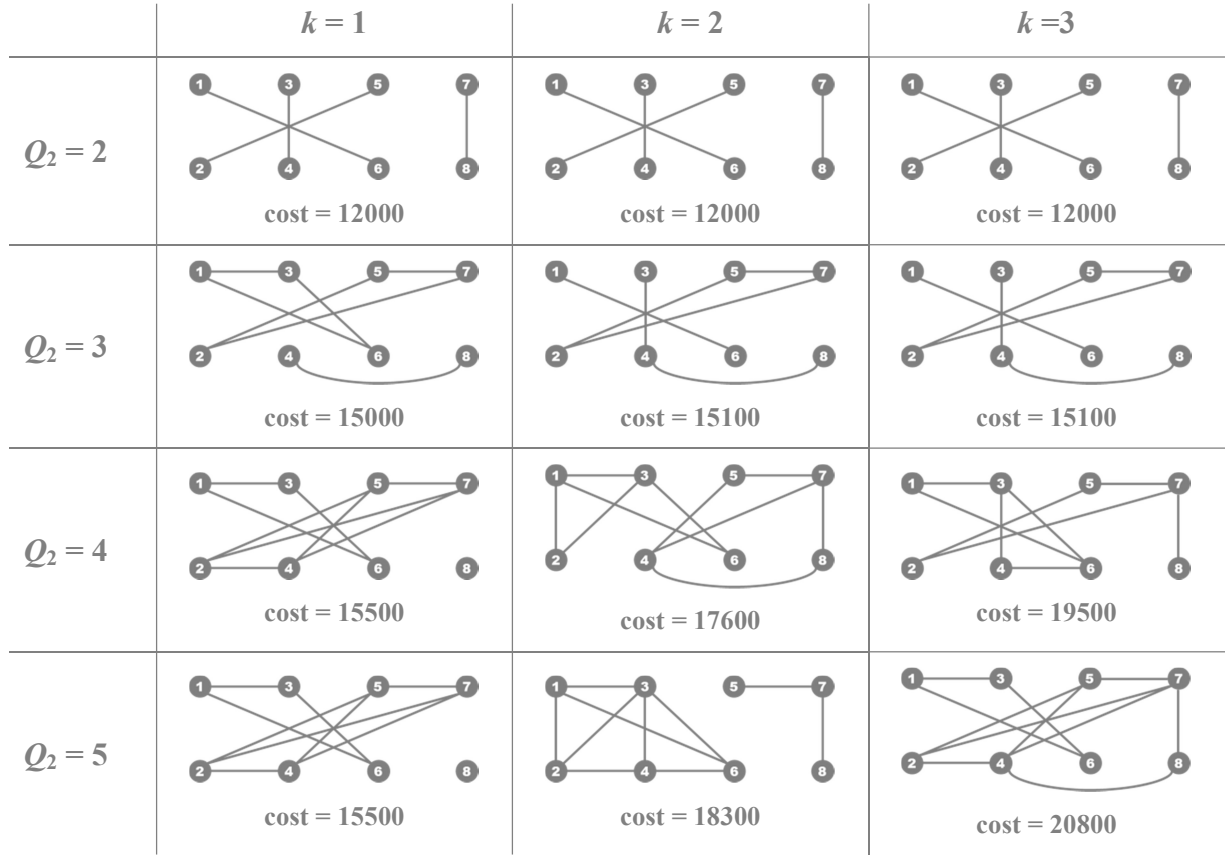

**Figure 4:** Optimum solutions of the Max-E$k$PP problem for $k = 1, 2, 3$ and $Q_2 = 2, 3, 4, 5$.

In this case, the solutions of Max-E$k$PP for $Q_2 \geq 6$ are the same as those obtained for $Q_2 = 5$. We could also observe the solutions involving the lower limit for the capacity in each component, or the upper limit for the number of components. However, the use of any of these additional topological constraints should be problem dependent, considering the practical application in hands. In addition, all these boundary conditions can be handled together, as long as the feasibility set is not empty.

## 3. Modeling the Max-E$k$PP problem

In the present section we propose mathematical formulations for the Max-E$k$PP problem. We start addressing the Max-E$k$PP for $k = 1$, that is, the Max-ECP problem. Then, we discuss the versions for $k > 1$. We also include a section for introducing the additional topological constraints.

## 3.1 The Max-E*k*PP for *k* = 1

As mentioned above, the Max-E*k*PP for *k* = 1 (or Max-ECP) has long been discussed in the literature. However, in all those cases, it has been modeled assuming that *G* is complete. The most popular formulation for the Max-ECP involves a single set of binary edge variables,

$$x_{ij} = \begin{cases} 1 & \text{if nodes } i \text{ and } j \text{ belong to the same component} \\ 0 & \text{otherwise} \end{cases}, \text{ for all } (i,j) \in E$$

and the model is (see, Grötschel and Wakabayashi 1989, 1990)

$$\max \sum_{(i,j) \in E} w_{ij} x_{ij}$$

(F1c)

s.t.
$$x_{ij} + x_{jk} - x_{ik} \leq 1, \text{ for all } (i,j),(j,k),(i,k) \in E \text{ with } i < j < k \quad (1)$$
$$x_{ij} - x_{jk} + x_{ik} \leq 1, \text{ for all } (i,j),(j,k),(i,k) \in E \text{ with } i < j < k \quad (2)$$
$$-x_{ij} + x_{jk} + x_{ik} \leq 1, \text{ for all } (i,j),(j,k),(i,k) \in E \text{ with } i < j < k \quad (3)$$
$$x_{ij} \in \{0,1\}, \text{ for all } (i,j) \in E \quad (4)$$

We denote this formulation by F1c (considering that *k* = 1 and *G* is complete). The set of constraints (1-3) is denoted by triangle inequalities, and they guarantee the transitivity propriety, establishing that if nodes *i* and *j* are in the same component (edge (*i,j*) is in the solution) and *j* and *k* also belong to the same component (edge (*j,k*) is in the solution), then the three nodes integrate the same component and edge (*i,k*) is forced to belong to the solution. These inequalities guarantee that each component represents a clique. In fact, due to the inclusionwise property, every subset of a clique is still a clique. So, a clique can be seen as a composition of triangles.

An alternative formulation is also proposed in Johnson et al. (1993). It uses edge/component and node/component variables, so requiring the prior definition of a set of components. The number of variables is much larger than in F1c, although requiring a smaller number of constraints.

Another alternative formulation is discussed in Mehrotra and Trick (1998). In this case, it considers an additional set of node/component variables and is based on the minimization of the sum of the weights of the edges that have end nodes in different components. However, this alternative formulation involves an even larger set of constraints.

In spite of the large number of inequalities involved in F1c, the model is still considered in recent works on the Max-ECP problem (see, e.g., Jaehn and Pesch 2013; Bettinelli et al.

2015; Brimberg, et al. 2015). Besides, as shown in Grötschel and Wakabayashi (1990), the linear programming relaxation (LP-relaxation) of F1c is "a quite reasonable" relaxation of the convex hull of the set of incidence vectors that characterize clique partitionings of $G$ (assumed complete). So, we follow the same motivation and conduct all the present work over the modeling structure of F1c.

If graph $G$ is sparse, the number of triangles decreases and the model should reflect it. The inclusion of the missing edges with sufficiently negative weights on the edges, proposed in Brimberg et al. (2015), does not allow profiting from the sparseness of the graph. On the other hand, if we remove from F1c all the variables associated to absent edges, the model provides a correct characterization of the Max-ECP problem. Yet, it includes a large number of redundant triangle constraints, namely all those involving 3 or 2 missing edges. We should also remove all the triangle inequalities involving a single missing edge and such that the missing edge does not correspond to the variable with negative sign. Thus, if we remove those constraints, we can obtain a much more compact formulation for the Max-ECP, depending on the sparseness of the graph. It uses that same set of $\{x_{ij}\}$ variables. We denote it by F1s.

(F1s)
$$\max \sum_{(i,j)\in E} w_{ij} x_{ij}$$
$$\text{s.t.} \quad x_{ij} + x_{jk} - x_{ik} \leq 1, \quad \text{for all } (i,j),(j,k),(i,k)\in E \text{ with } i<j<k \quad (5)$$
$$x_{ij} - x_{jk} + x_{ik} \leq 1, \quad \text{for all } (i,j),(j,k),(i,k)\in E \text{ with } i<j<k \quad (6)$$
$$-x_{ij} + x_{jk} + x_{ik} \leq 1, \quad \text{for all } (i,j),(j,k),(i,k)\in E \text{ with } i<j<k \quad (7)$$
$$x_{ik} + x_{jk} \leq 1, \quad \text{for all } (i,k),(j,k)\in E \text{ and } (i,j)\notin E \text{ with } i<j<k \quad (8)$$
$$x_{ij} \in \{0,1\}, \quad \text{for all } (i,j)\in E \quad (9)$$

The set of constraints (5-8) are obtained from the triangle inequalities (1-3) after removing the previously mentioned redundant constraints. If we ignore inequalities (8) from F1s, the feasibility set includes other solutions besides clique partitionings. For instance, using again the example in Figure 1, the solution returned by F1s without constraints (8) is the entire graph $G$.

## 3.2 The Max-E*k*PP for *k* > 1

For modeling the Max-E*k*PP with $k > 1$, we introduce an additional set of variables characterizing the missing edges in $G$. The new variables are defined by

$$v_{ij} = \begin{cases} 1 & \text{if nodes } i \text{ and } j \text{ belong to the same component} \\ 0 & \text{otherwise} \end{cases} \text{, for all } (i,j) \in E^c$$

We further use the following notation

$$\{x_{ij}, v_{ij}\} = \begin{cases} x_{ij} & \text{if } (i,j) \in E \\ v_{ij} & \text{if } (i,j) \in E^c \end{cases} \text{, for all } (i,j) \in E \cup E^c$$

Hence, the formulation for the Max-E*k*PP with $k > 1$, denoted by F*k*s, is

max $\sum_{(i,j) \in E} w_{ij} x_{ij}$

s.t.  $\{x_{ij}, v_{ij}\} + \{x_{jk}, v_{jk}\} - \{x_{ik}, v_{ik}\} \leq 1$ , for all $(i,j),(j,k),(i,k) \in E \cup E^c$ with $i < j < k$ (10)

(F*k*s) $\{x_{ij}, v_{ij}\} - \{x_{jk}, v_{jk}\} + \{x_{ik}, v_{ik}\} \leq 1$ , for all $(i,j),(j,k),(i,k) \in E \cup E^c$ with $i < j < k$ (11)

$-\{x_{ij}, v_{ij}\} + \{x_{jk}, v_{jk}\} + \{x_{ik}, v_{ik}\} \leq 1$ , for all $(i,j),(j,k),(i,k) \in E \cup E^c$ with $i < j < k$ (12)

$\sum_{(i,j) \in \delta^c(i)} v_{ij} \leq k - 1$ , for all $i \in V$ (13)

$x_{ij} \in \{0,1\}$ , for all $(i,j) \in E$ (14)

$v_{ij} \in \{0,1\}$ , for all $(i,j) \in E^c$ (15)

Using the additional set of variables $\{v_{ij}\}$, we can characterize each component as a clique, composed by edges from $E$ and from $E^c$. Thus, we can resort to the structure of the triangle inequalities (1-3) for characterizing partitions into cliques on a complete graph. This way, as we can control the number of missing edges belonging to the solution (variables $\{v_{ij}\}$), we can put an upper limit on the number missing edges incident to a given node $i \in V$, imposed by inequalities (13).

For the particular case with $k = 2$, the model has many redundant constraints. In fact, due to constraints (13), all the triangle inequalities on $(i,j,k) \in V^3$ such that the three edges belong to set $E^c$ are redundant. The same way, and also due to (13), among those involving two missing edges, all the triangle inequalities $v_{ij} + v_{jk} - x_{ik} \leq 1$ for $(i,j),(j,k) \in E^c$ and $(i,k) \in E$ are redundant, as well.

The number of variables in F$k$s is $(n \cdot (n-1))/2$ and the number of constraints is $3 \cdot \binom{n}{3} + n$, when $k \geq 3$. This is almost the same size of model F1c for the Max-ECP problem on complete graphs. As mentioned above, the model has a smaller number of inequalities when $k = 2$.

Model F$k$s could have been defined just using variables $x$ notation, extending the original set to the edges in the complementary graph $G^c$. This alternative characterization could possibly simplify the notation. However, we opted to use variables $v$ in order to clearly distinguish the two sets of edges, trying to simplify the exposition.

### 3.3 Additional topological constraints

From a non-hierarchical clustering perspective, and assuming that the topological structure matters, finding a partition of $G$ into cliques or into $k$-plexes may also come with additional constraints. The most usually discussed constraints within the Max-ECP problem involve upper and lower bounds on the sum of the weights of the nodes in each component (see, Johnson et al. 1993; Ferreira et al. 1996; Mehrotra and Trick 1998; Ji and Mitchell 2007; Oosten et al. 2007); and/or an upper limit on the number of components in the solution (see, Mehrotra and Trick 1998).

We can model lower and upper bound limits on the sum of the weights of the nodes in each component using the set of variables $\{x_{ij}\}$ and $\{v_{ij}\}$ in the previously described models. Considering $lb$ as the lower bound and $ub$ the upper bound, the mentioned constraints are,

- for models F1c and F1s:

$$\text{lower bound constraints} \quad \sum_{j \in \delta(i)} q_j x_{ij} \geq lb - q_i \quad , \quad \text{for all } i \in V \quad (16)$$

$$\text{upper bound constraints} \quad \sum_{j \in \delta(i)} q_j x_{ij} \leq ub - q_i \quad , \quad \text{for all } i \in V \quad (17)$$

- for models F$k$s:

$$\text{lower bound constraints} \quad \sum_{j \in V \setminus \{i\}} q_j \{x_{ij}, v_{ij}\} \geq lb - q_i \quad , \quad \text{for all } i \in V \quad (18)$$

$$\text{upper bound constraints} \quad \sum_{j \in V \setminus \{i\}} q_j \{x_{ij}, v_{ij}\} \leq ub - q_i \quad , \quad \text{for all } i \in V \quad (19)$$

Within the Max-ECP problem literature, lower bounding constraints similar to those in (16) were considered in Oosten et al. (2007), while alternative versions of the upper bounding constraints (17) were discussed in Ferreira et al. (1996) and Ji and Mitchell (2007) using node variables.

On the other hand, the set of variables considered in formulation F1c does not allow an easy adaptation for modeling a limit on the number of components, as mentioned in Mehrotra and Trick (1998). So, we consider an additional set of node/component variables that relate each node to each of the available components, that is

$$z_i^p = \begin{cases} 1 & \text{if node } i \text{ belongs to component } p \\ 0 & \text{otherwise} \end{cases}, \text{ for all } i \in V \text{ and } p \in \{1,\ldots,P\}$$

with $\{1,\ldots,P\}$ the set of components and $P$ the upper bound for the number of components. Thus, the following constraints model the intended limitation

- for models F1c and F1s:

$$\sum_{p=1}^{P} z_i^p = 1 \quad , \quad \text{for all } i \in V \tag{20}$$

$$z_i^p + z_j^p \leq 1 + x_{ij} \quad , \quad \text{for all } (i,j) \in E \text{ and } p \in \{1,\ldots,P\} \tag{21}$$

$$z_i^p + z_j^p \leq 1 \quad , \quad \text{for all } (i,j) \in E^c \text{ and } p \in \{1,\ldots,P\} \tag{22}$$

- for models F$k$s:

$$\sum_{p=1}^{P} z_i^p = 1 \quad , \quad \text{for all } i \in V \tag{23}$$

$$z_i^p + z_j^p \leq 1 + \{x_{ij}, v_{ij}\} \quad , \quad \text{for all } i,j \in V \text{ with } i<j \text{ and } p \in \{1,\ldots,P\} \tag{24}$$

Constraints (20) guarantee that each node belongs to one of the components. Inequalities (21) impose that two nodes can only be in the same component if there is an edge linking them in $G$; otherwise, when the edge does not exists, the two nodes cannot share a component, as stated in (22).

Then, for those addressing models F$k$s, constraints (23) are the same as (20); while inequalities (24) integrate (21) and (22), stating that two nodes may belong to the same component if they are linked by and edge or if they are related by a missing edge selected by the solution.

As mentioned above, upper limit constraints on the number of components were also discussed in Ferreira et al. (1996) and Mehrotra and Trick (1998). They consider the same set of variables $\{z_i^p\}$, using the range of variation of the extra $p$ index to set the limit on the number of components in the solution.

# 4. Computational tests

In the present section we discuss the models proposed in Section 3. This discussion is conducted in two parts:

*i*) exploring the applicability of the problem, including some of the additional topological constraints introduced in Subsection 3.3; and

*ii*) analyzing the performance of the models proposed in Subsections 3.1 and 3.2 using a commercial solver.

These two parts are discussed in the forthcoming subsections.

The models were solved using ILOG/CPLEX 11.2 and all experiments were performed under Microsoft Windows 7 operating system on an Intel Core i7-2600 with 3.40 GHz and 8 GB RAM. When running the mixed integer programming (MIP) algorithm of CPLEX we used most default settings, which involve an automatic procedure that uses the best rule for variable selection and the best-bound search strategy for node selection in the branch-and-bound tree. We have set an upper time limit of 10800 seconds for each test. The times are reported in seconds.

## 4.1 Applying the Max-ECP and the Max-E*k*PP problems

In order to discuss the applicability of problems Max-ECP and Max-E*k*PP, we consider two types of biological networks involving metabolic reactions' interactions and metabolite's interactions. These are *Saccharomyces cerevisiae* metabolic networks, taken from Förster et al. (2003). The original data involves 1393 metabolic reactions (we substituted reciprocal (bidirectional) reactions into two single-direction reactions) that use 991 metabolites. Each metabolic reaction is a chemical pathway that uses reactants to generate products. Both reactants and products are metabolites being shared among reactions. For instance, we can characterize a reaction as: A + B → C + D, meaning that metabolites A and B are reactants producing metabolites C and D (as products). The same way, we can have another reaction characterized by: A + C + E → B + F. The two reactions share three different metabolites: A, B and C. Using the entire set of metabolic reactions from the *Saccharomyces cerevisiae* data, described in the mention paper, we built the following network structures:

- <u>SC-NIP-m-t*r*</u> (for $r = 1,\ldots,5$): each node in $V$ represents a metabolite (excluding isolated ones) and each edge in $E$ represents a pair of metabolites that share at least $r$ reaction; thus ($i,j$)∈$E$ if there are at least $r$ common reactions that include both metabolites $i$ and $j$ (no matter the side: reactant or product). Weight $c_{ij}$ represents the number of reactions sharing the two metabolites $i$ and $j$.

- <u>SC-NIP-r-t$m$</u> (for $m = 1,…,5$): each node in $V$ represents a reaction (excluding isolated ones) and each edge in $E$ represents a pair of reactions sharing at least $m$ metabolites; thus ($i,j$)∈$E$ if there are at least $m$ common metabolites among the two reactions $i$ and $j$ (no matter the side the metabolites appear in). Weight $c_{ij}$ represents the number of metabolites sharing the two reactions $i$ and $j$.

Isolated nodes were removed from the original data, in all instances under discussion. Table 1 indicates the number of nodes, the number of edges and the density of each of the proposed graphs. $d$ denotes the density of the graph, with $d = 2|E|/(n \cdot (n-1))$.

| Instance | $n$ | $|E|$ | $d$ | Instance | $n$ | $|E|$ | $d$ |
|---|---|---|---|---|---|---|---|
| SC-NIP-m-t1 | 991 | 4161 | 0.00848 | SC-NIP-r-t1 | 1393 | 56276 | 0.05804 |
| SC-NIP-m-t2 | 602 | 1520 | 0.00840 | SC-NIP-r-t2 | 1183 | 17776 | 0.02542 |
| SC-NIP-m-t3 | 177 | 269 | 0.01727 | SC-NIP-r-t3 | 663 | 1782 | 0.00812 |
| SC-NIP-m-t4 | 129 | 166 | 0.02011 | SC-NIP-r-t4 | 377 | 321 | 0.00453 |
| SC-NIP-m-t5 | 75 | 84 | 0.03027 | SC-NIP-r-t5 | 45 | 27 | 0.02727 |

**Table 1:** Characterization of the two types of *Saccharomyces cerevisiae* graphs under discussion.

We have applied problem Max-E$k$PP for all the reported instances, considering $k = 1, 2, 3$. When $k = 1$, the problem is the Max-ECP. Thus, we have used model F$k$s with $k = 1, 2, 3$ for solving the problems and resorted to CPLEX to solve the models.

Tables 2 and 3 show the results for the classes of instances SC-NIP-m and SC-NIP-r, respectively; concerning problems Max-ECP and Max-E$k$PP for $k = 2, 3$. They report:

- the optimum values or the best feasible solutions (opt/best) found by the branch-and-bound;
- the duality gap at termination, defined by d_gap = ((UB – LB) / UB)*100, in percent, with UB indicating the best upper bound and LB indicating the best lower bound or the optimum, found by the branch-and-bound;
- "time" reports the branch-and-bound execution time (in seconds);
- "comp" indicates the number of components in the solution;
- "largest" is the number of nodes in the largest cardinality component in the solution;
- "singlt" indicates the proportion of singletons over the total number of nodes (in percent).

When the optimum is not attained within the given time limit, the d_gap is not null and the time is denoted by " > ". When memory is insufficient for building the model or insufficient for reading it, we put the following indication: "o.m." (out-of memory).

| k | Instance | opt/best | d_gap | time | comp | largest | singlt |
|---|---|---|---|---|---|---|---|
| 1 | SC-NIP-m-t1 | 1866 | 0.00 | 2296.94 | 340 | 8 | 16.65 |
|   | SC-NIP-m-t2 | 1538 | 0.00 | 1.25 | 206 | 7 | 19.77 |
|   | SC-NIP-m-t3 | 910 | 0.00 | 0.02 | 48 | 5 | 35.03 |
|   | SC-NIP-m-t4 | 831 | 0.00 | 0.00 | 38 | 5 | 31.78 |
|   | SC-NIP-m-t5 | 723 | 0.00 | 0.00 | 27 | 4 | 21.33 |
| 2 | SC-NIP-m-t1 | o.m. | ---- | ---- | ---- | ---- | ---- |
|   | SC-NIP-m-t2 | o.m. | ---- | ---- | ---- | ---- | ---- |
|   | SC-NIP-m-t3 | 1021 | 0.00 | 50.43 | 47 | 7 | 28.81 |
|   | SC-NIP-m-t4 | 907 | 0.00 | 3.03 | 39 | 5 | 24.03 |
|   | SC-NIP-m-t5 | 801 | 0.00 | 0.20 | 26 | 4 | 17.33 |
| 3 | SC-NIP-m-t1 | o.m. | ---- | ---- | ---- | ---- | ---- |
|   | SC-NIP-m-t2 | o.m. | ---- | ---- | ---- | ---- | ---- |
|   | SC-NIP-m-t3 | o.m. | ---- | ---- | ---- | ---- | ---- |
|   | SC-NIP-m-t4 | o.m. | ---- | ---- | ---- | ---- | ---- |
|   | SC-NIP-m-t5 | 887 | 0.00 | 34.20 | 24 | 5 | 16.00 |

**Table 2:** Optimum solutions for the class of instances SC-NIP-m, using models F$k$s, $k$=1, 2, 3.

| k | Instance | opt/best | d_gap | time | comp | largest | singlt |
|---|---|---|---|---|---|---|---|
| 1 | SC-NIP-r-t1 | o.m. | ---- | ---- | ---- | ---- | ---- |
|   | SC-NIP-r-t2 | 34576 | 0.00 | 4.26 | 262 | 121 | 4.40 |
|   | SC-NIP-r-t3 | 5411 | 0.00 | 0.08 | 245 | 39 | 3.92 |
|   | SC-NIP-r-t4 | 1232 | 0.00 | 0.00 | 170 | 11 | 1.86 |
|   | SC-NIP-r-t5 | 140 | 0.00 | 0.02 | 21 | 4 | 2.22 |
| 2 | SC-NIP-r-t1 | o.m. | ---- | ---- | ---- | ---- | ---- |
|   | SC-NIP-r-t2 | o.m. | ---- | ---- | ---- | ---- | ---- |
|   | SC-NIP-r-t3 | 3183 | 71.62 | > | 176 | 35 | 31.37 |
|   | SC-NIP-r-t4 | 1245 | 0.00 | 6.40 | 170 | 11 | 1.06 |
|   | SC-NIP-r-t5 | 140 | 0.00 | 0.01 | 21 | 4 | 2.22 |
| 3 | SC-NIP-r-t1 | o.m. | ---- | ---- | ---- | ---- | ---- |
|   | SC-NIP-r-t2 | o.m. | ---- | ---- | ---- | ---- | ---- |
|   | SC-NIP-r-t3 | o.m. | ---- | ---- | ---- | ---- | ---- |
|   | SC-NIP-r-t4 | o.m. | ---- | ---- | ---- | ---- | ---- |
|   | SC-NIP-r-t5 | 140 | 0.00 | 0.14 | 21 | 4 | 2.22 |

**Table 3:** Optimum/best solutions for the class of instances SC-NIP-r, using models F$k$s, $k$=1, 2, 3.

From the computational stand point, all instances have been solved to optimality for the version with $k = 1$ (Max-ECP), except the largest example involving metabolic reactions' interactions (SC-NIP-r-t1). For this version the model is more compact. In effect, model's F1s number of variables and constraints depends on the sparsity of the graph. The version with $k = 2$ is harder as it uses model F2s that includes a larger number of variables and constraints. In this case, we have not been able to read the models involving the two largest instances in both classes (reactions' interactions and metabolites' interactions); and instance SC-NIP-r-t3 could not reach the optimum within the given time limit. The harder task, however, was observed solving the version with $k = 3$, which requires model F3s. In this case, we could only solve the smaller instances SC-NIP-m-t5 and SC-NIP-r-t5. In fact, the number of variables and constraints in model F3s (and for all models F$k$s with $k \geq 3$) is almost the same as if we were dealing with the cliques partitioning problem on the complete graph, so almost the same size as model F1c. So, it is not a surprise to stop sooner when trying to solve increasingly larger sized instances.

When looking to the solutions, and starting with the networks involving metabolite's interactions (instances SC-NIP-m), the heaviest components reveal some of the most shared metabolites. If ignoring the singletons and the low heavy components, the solutions show a modular structure of the metabolites in the system. To further stress this observation, we report in Table 4 the main information on the three heaviest components in each of the optimum solutions reported in Table 2. The list of the metabolites involved in each of the mentioned components is shown in Table A1 in the Appendix.

| $k$ | Instance | heaviest components | total weight | total different reactions | cardinality |
|---|---|---|---|---|---|
| 1 | SC-NIP-m-t1 | 1st | 343 | 191 | 8 |
|   |             | 2nd | 102 | 90  | 4 |
|   |             | 3rd | 87  | 63  | 6 |
|   | SC-NIP-m-t2 | 1st | 320 | 186 | 7 |
|   |             | 2nd | 102 | 90  | 4 |
|   |             | 3rd | 76  | 62  | 5 |
|   | SC-NIP-m-t3 | 1st | 242 | 136 | 5 |
|   |             | 2nd | 96  | 90  | 3 |
|   |             | 3rd | 56  | 54  | 3 |
|   | SC-NIP-m-t4 | 1st | 242 | 136 | 5 |
|   |             | 2nd | 90  | 90  | 2 |

|   |           | rank | total weight | total different reactions | cardinality |
|---|-----------|------|--------------|---------------------------|-------------|
|   |           | 3rd  | 50           | 16                        | 4           |
|   | SC-NIP-m-t5 | 1st | 214          | 125                       | 4           |
|   |           | 2nd  | 90           | 90                        | 2           |
|   |           | 3rd  | 45           | 45                        | 2           |
| 2 | SC-NIP-m-t3 | 1st | 342          | 187                       | 7           |
|   |           | 2nd  | 96           | 90                        | 3           |
|   |           | 3rd  | 59           | 37                        | 4           |
|   | SC-NIP-m-t4 | 1st | 276          | 184                       | 5           |
|   |           | 2nd  | 90           | 90                        | 2           |
|   |           | 3rd  | 59           | 37                        | 4           |
|   | SC-NIP-m-t5 | 1st | 262          | 178                       | 4           |
|   |           | 2nd  | 90           | 90                        | 2           |
|   |           | 3rd  | 59           | 37                        | 4           |
| 3 | SC-NIP-m-t5 | 1st | 326          | 181                       | 5           |
|   |           | 2nd  | 126          | 116                       | 4           |
|   |           | 3rd  | 79           | 37                        | 5           |

**Table 4:** Three heaviest components' information of the optimum solutions reported in Table 2.

In Table 4 we distinguish the total number of reactions sharing any pair of metabolites in the component (column "total weight"), from the total number of distinct reactions involved in the component (column "total different reactions"). Column "cardinality" indicates the number of metabolites in the component.

The selected components in Table 4 show that there are a very small number of metabolites acting together in a large number of reactions. Some of these reactions share more than a pair of metabolites in the same component. However, the total number of distinct reactions is still large compared to the number of metabolites involved.

To exemplify, Figures 5, 6 and 7 show the mentioned heaviest components involving instance SC-NIP-m-t5 for $k = 1$, 2 and 3, respectively. We recall that each edge in the graph associated to instance SC-NIP-m-t5 links a pair of metabolites sharing at least 5 metabolic reactions.

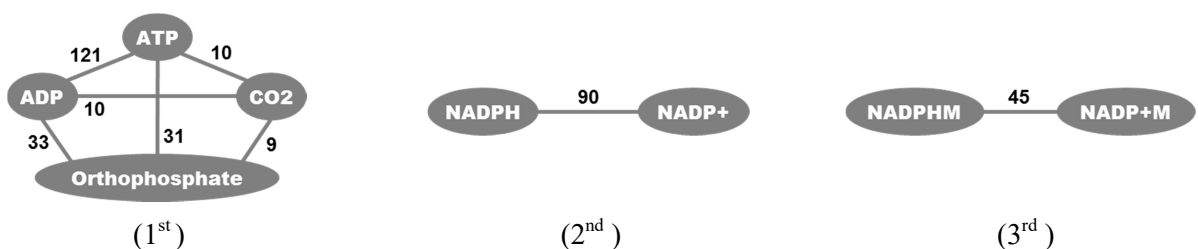

**Figure 5:** Three heaviest components in the optimum solution of instance SC-NIP-m-t5 for $k = 1$.

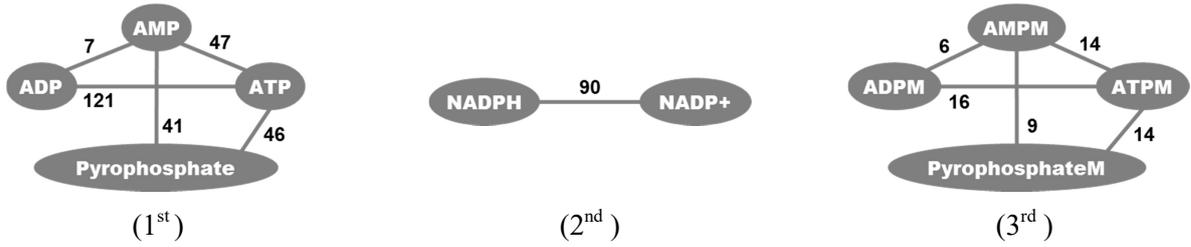

**Figure 6:** Three heaviest components in the optimum solution of instance SC-NIP-m-t5 for $k = 2$.

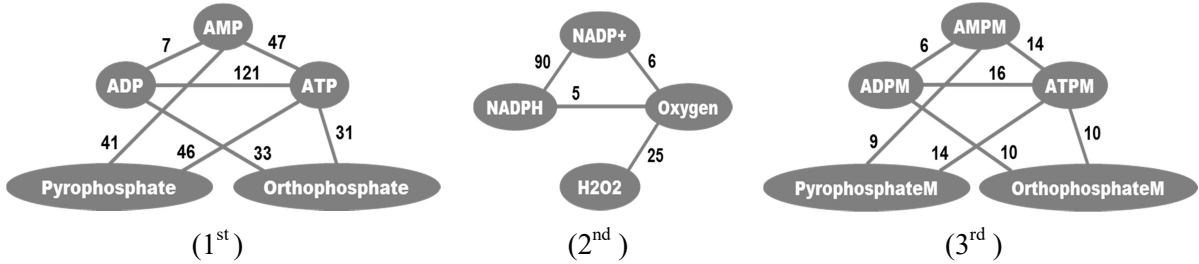

**Figure 7:** Three heaviest components in the optimum solution of instance SC-NIP-m-t5 for $k = 3$.

The exemplified solutions show, in general, that the components get larger when the $k$-plexes become more relaxed, bringing additional metabolites into the modules. In addition and when considering the heaviest components, the metabolites involved can be seen as hubs in the entire system of metabolic reactions.

On the other hand, the solutions on the networks involving metabolic reactions' interactions (instances SC-NIP-r) present a partition of the entire set of reactions into modules. This modular structure of the system may show a higher level picture of the system, involving a modular structure on macro-pathways, namely on the heaviest components. In this case, the heaviest components are quite large, compared with the former approach, as indicated in Table 3 and shown in Table 5 for the three heaviest components in each optimum solution found using the instances graphs SC-NIP-r. We have omitted the feasible solution found for instance SC-NIP-r-t3 with $k = 2$ due to its large gap.

| $k$ | Instance | heaviest components | total weight | total different metabolites | cardinality |
|---|---|---|---|---|---|
| 1 | SC-NIP-r-t2 | 1st | 15188 | 68 | 121 |
|   |             | 2nd | 8077  | 56 | 90  |
|   |             | 3rd | 2756  | 13 | 45  |
|   | SC-NIP-r-t3 | 1st | 2296  | 13 | 39  |
|   |             | 2nd | 1399  | 26 | 30  |

|       |            | 3rd | 122 | 12 | 9  |
|-------|------------|-----|-----|----|----|
| 1, 2  | SC-NIP-r-t4 | 1st | 220 | 4  | 11 |
|       |            | 2nd | 150 | 10 | 9  |
|       |            | 3rd | 89  | 7  | 7  |
| 1, 2, 3 | SC-NIP-r-t5 | 1st | 30  | 5  | 4  |
|       |            | 2nd | 6   | 6  | 2  |
|       |            | 3rd | 6   | 6  | 2  |

**Table 5:** Three heaviest components' information of the optimum solutions reported in Table 3.

As mentioned above, the three heaviest components in the solutions reported in Table 5 involve a very large number of metabolic reactions, also sharing a large number of metabolites. Yet, the entire sum of metabolites shared among the reactions (characterizing the components' weight) involves many repetitions. In effect, when looking just for the set of distinct metabolites being shared in each component, its cardinality becomes much smaller.

Still within the networks involving metabolic reactions' interactions (instances SC-NIP-r), we have also explored the partitioning problem considering one of the types of topological constraints described in Subsection 3.3. So, considering the particular instance SC-NIP-r-t3 with $k = 1$, we have included the additional condition that sets an upper bound for the cardinality in each component. To this purpose, we took model F1s addressing problem Max-ECP with the additional set of constraints (17), and assuming that $q_i = 1$ for all $i \in V$. The problem was discussed considering the upper bound limit set to $ub = 10$ and 5, that is, forcing all components in any feasible partition to have no more than 10 and 5 metabolic reactions, respectively. The purpose is to obtain more balanced solutions concerning the cardinality of the components. Notice that the mentioned instance cannot be solved when imposing a lower limit equal to 2 or higher in the cardinality of the components. This is due to the sparsity of the associated graph.

Table 6 shows the main information concerning the best solutions found for the two cases, using the same notation considered in Table 3.

| $k$ | Instance    | $ub$ | opt/best | d_gap | time | comp | largest | singlt |
|-----|-------------|------|----------|-------|------|------|---------|--------|
| 1   | SC-NIP-r-t1 | 5    | 1991     | 1.10  | >    | 261  | 5       | 3.77   |
|     |             | 10   | 2805     | 0.46  | >    | 249  | 10      | 4.07   |

**Table 6:** Best solutions for instance SC-NIP-r-t3 with $k = 1$ and $ub = 5$ and 10.

We have not reached the optimum in both cases within the given time limit. However, the gaps are relatively small. As expected, the number of components increases when the upper bound on the components' cardinalities gets tighter. Figure 8 describes the number of components for the range of cardinalities in the solution of the unrestricted version and in the solutions of the two restricted cases, for $ub$ = 5 and 10.

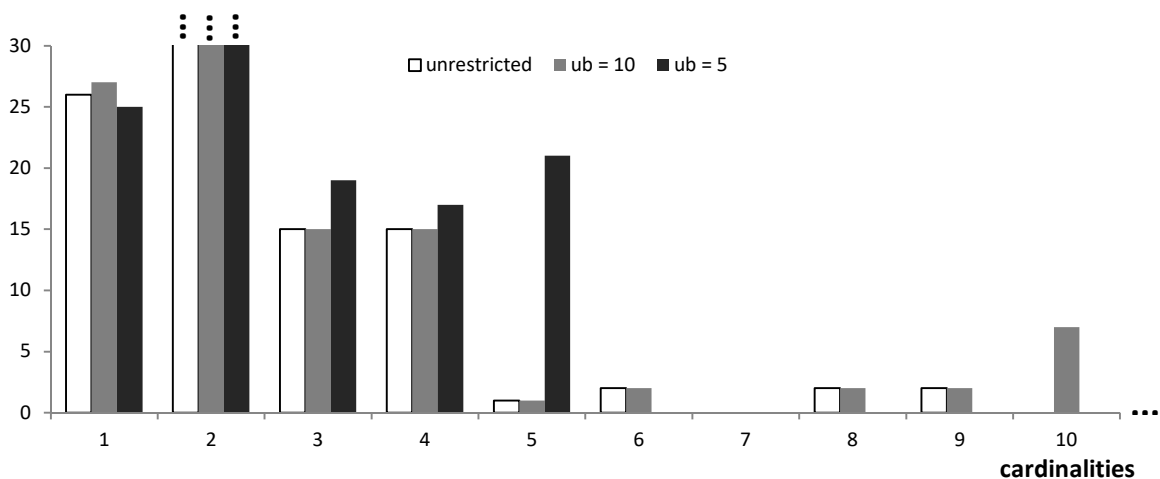

**Figure 8:** Number of components for the range of cardinalities in the solutions for the unrestricted version and for the cases with $ub$ = 5 and 10.

The horizontal axis is set only up to 10. The missing information concerns the unrestricted version, involving just one component with cardinality 30 and one component with cardinality 39. Also, the vertical axis varies only up to 30, cutting just the number of components of cardinality 2. In fact, the number of components of that type is very large, being equal to 206, 205 and 204 for the unrestricted version and for the restricted cases with $ub$ = 10 and 5, respectively. The number of nodes (metabolic reactions) involved in these components represents 62.14%, 61.84% and 61.54% of the entire number of nodes in the graph, respectively. As observed in the graphic in Figure 8, apart from those with cardinality 2, the number of components in each of the other cardinalities gets more balanced as the upper bound becomes smaller. Actually, in the case with $ub$ = 5, the number of components with cardinalities 3, 4 and 5 is equal to 19, 17 and 21, respectively.

## 4.2 Computational performance of models F1s and Fks

This subsection intends to perform some computational test on the models proposed in Subsections 3.1 and 3.2 for sparse graphs. To this purpose, we consider a selection of benchmark instances taken from the DIMACS database. We took all instances with $|V| \leq 100$ and all those with $100 < |V| \leq 200$ and density ($d$) at most 0.25. These instances belong to the c-fat, MANN, hamming and johnson families. The original DIMACS instances do not include weights on edges. So, we followed the weighting strategy proposed in Pullan (2008), setting $w_{ij} = ((i + j) \mod 200) + 1$. We do not report the other DIMACS instances' results because we were not able to answer them.

Table 7 presents the computational results' information using models F$k$s for $k = 1, 2, 3$. The columns' labels are the same as those considered in Tables 2 and 3. The notation "-----", in this case, means that we have not found a single feasible solution.

| Instance | $n$ | $d$ | $k$ | opt/best | d_gap | time | comp | largest | singlt |
|---|---|---|---|---|---|---|---|---|---|
| c-fat200-1 | 200 | 0.077 | 1 | 98711 | 0.00 | 47.08 | 19 | 12 | 0.00 |
| | | | 2 | 98711 | 0.00 | 567.44 | 19 | 12 | 0.00 |
| | | | 3 | ----- | ----- | ----- | ----- | ----- | ----- |
| c-fat200-2 | 200 | 0.163 | 1 | 213248 | 0.00 | 0.22 | 9 | 24 | 0.00 |
| | | | 2 | 213248 | 0.00 | 47.28 | 9 | 24 | 0.00 |
| | | | 3 | ----- | ----- | ----- | ----- | ----- | ----- |
| hamming-6-2 | 64 | 0.905 | 1 | 65472 | 0.00 | 0.20 | 2 | 32 | 0.00 |
| | | | 2 | 65472 | 6.25 | > | 2 | 32 | 0.00 |
| | | | 3 | 65472 | 23.30 | > | 2 | 32 | 0.00 |
| hamming-6-4 | 64 | 0.349 | 1 | 6336 | 0.00 | 0.34 | 16 | 4 | 0.00 |
| | | | 2 | 6966 | 149.85 | > | 13 | 6 | 1.56 |
| | | | 3 | 4567 | 307.67 | > | 21 | 4 | 0.00 |
| johnson8-2-4 | 28 | 0.556 | 1 | 1260 | 0.00 | 0.06 | 7 | 4 | 0.00 |
| | | | 2 | 1355 | 57.89 | > | 6 | 5 | 0.00 |
| | | | 3 | 1996 | 36.67 | > | 4 | 8 | 0.00 |
| johnson8-4-4 | 70 | 0.768 | 1 | 27864 | 17.57 | > | 6 | 14 | 0.00 |
| | | | 2 | 12770 | 438.28 | > | 11 | 9 | 2.86 |
| | | | 3 | 12948 | 463.40 | > | 12 | 8 | 0.00 |
| MANN_a9 | 45 | 0.927 | 1 | 14868 | 0.00 | 1215.34 | 3 | 16 | 0.00 |
| | | | 2 | 23047 | 1.82 | > | 2 | 26 | 0.00 |
| | | | 3 | 33660 | 0.00 | 319.24 | 2 | 36 | 0.00 |

**Table 7:** Computational tests' information on the selected DIMACS instances using models F$k$s, $k = 1, 2, 3$.

All instances under discussion were solved to optimality when addressing the Max-ECP problem, except for instance johnson8-4-4. This is the larger dimensional example with density $d > 0.25$.

For the other problems, with $k \geq 2$, addressing the Max-E$k$PP, things become much harder, especially when $k$ gets larger, in general. The main exception is instance MANN_a9. A reason for this exception can be related with its very high density. This observation seems contradictory with another one in the previous paragraph. However, when the graph is very dense, larger values of $k$ allow higher relaxed components, benefiting an easier coverage of the entire set of nodes using a smaller number of components. This is probably the same reason for justifying the smaller gaps observed for instance hamming-6-2 (denser) when compared with those obtained for instance hammin-6-4 (sparser). A similar observation applies between instances c-fat200-2 (denser) and c-fat200-1 (sparser).

From an empirical stand point, these tests suggest that solving the Max-E$k$PP is harder for sparser graphs. They also suggest that the versions with $k > 1$ are probably much harder than the Max-ECP, worsening when $k$ increases.

## 5. Conclusions

This paper discusses the maximum edge weight $k$-plex partitioning problem on sparse graphs. When $k = 1$, it becomes the well-known maximum edge weight clique partitioning.

To our best knowledge, the versions with $k \geq 2$ have never been discussed in the literature. So, we present the first formulations for those versions of the problem. In the computational tests performed, all the models were solved using a commercial branch-and-bound based solver (IBM/CPLEX) with no further techniques.

Although short, this computational experience was sufficient to confirm our worst expectations, showing that the Max-E$k$PP problem is probably much harder, from an empirical stand point, than the Max-ECP, and worsening when $k$ increases.

These tests involved instances with up to 200 nodes, being solved to optimality, in general, when addressing the Max-ECP problem. Yet, and for the same instances, we seldom solved to optimality the Max-E$k$PP within the given time limit. In fact, in most cases, we were probably far from the optimums.

An immediate conclusion is that the solving techniques require more sophisticated methodologies, namely resorting to decomposition methods or cutting plane processes, starting from relaxed versions of the original models.

Other relevant contributions may involve alternative ideas for modeling the Max-E$k$PP; or the identification of additional non-trivial cuts to further strengthen the models' linear programming relaxation polyhedron.

In addition, considering the potential practical relevancy of the problem, and reminding that many practical applications involve large sized graphs, we also detach the importance of building heuristic algorithms for the Max-E$k$PP.

Nevertheless, even considering the proposed models' limitations, we used them to discuss a real-world applied problem on two biological networks involving metabolic reactions' interactions and metabolite's interactions, respectively. Although very sparse, these graphs were much larger than those used in the previously mentioned computational tests. These examples were used to detach the applicability of the problems.

We should not ignore, however, that having a formulation and a solver engine can be an easier way for many researchers to first apply the Max-E$k$PP problem.

## Acknowledgments

This work was partially supported by the Portuguese National Funding: Fundação para a Ciência e a Tecnologia – FCT (project UID/MAT/04561/2013).

## Appendix

Table A1 presents the list of the metabolites involved in each of the components reported in Table 4, addressing the class of instances on networks involving metabolite's interactions (instances SC-NIP-m).

|   |   | heaviest components | | |
|---|---|---|---|---|
| $k$ | Instance | 1st | 2nd | 3rd |
| 1 | SC-NIP-m-t1 | ADP<br>AMP<br>ATP<br>L=Aspartate<br>L=Glutamate<br>L=Glutamine<br>NH3<br>Orthophosphate | L=2=Aminoadipate_6=semialdehyde<br>N6=(L=1,3=Dicarboxypropyl)=L=lysine<br>NADP+<br>NADPH | Acetyl=CoA<br>Acyl=carrier_protein<br>CO2<br>CoA<br>Hexadecanoyl=CoA<br>Malonyl=CoA |
|   | SC-NIP-m-t2 | ADP<br>AMP<br>ATP<br>L=Glutamate<br>L=Glutamine<br>NH3<br>Orthophosphate | L=2=Aminoadipate_6=semialdehyde<br>N6=(L=1,3=Dicarboxypropyl)=L=lysine<br>NADP+<br>NADPH | Acetyl=CoA<br>Acyl=carrier_protein<br>CO2<br>CoA<br>Hexadecanoyl=CoA<br>Malonyl=CoA |
|   | SC-NIP-m-t3 | ADP<br>ATP<br>L=Glutamate<br>L=Glutamine<br>Orthophosphate | L=2=Aminoadipate_6=semialdehyde<br>NADP+<br>NADPH | 5=Phospho=alpha=D=ribose_1=diphosphate<br>AMP<br>Pyrophosphate |
|   | SC-NIP-m-t4 | ADP<br>ATP<br>L=Glutamate<br>L=Glutamine<br>Orthophosphate | NADP+<br>NADPH | ADPM<br>ATPM<br>CoAM<br>OrthophosphateM |
|   | SC-NIP-m-t5 | ADP<br>ATP<br>CO2<br>Orthophosphate | NADP+<br>NADPH | NADPHM<br>NADP+M |
| 2 | SC-NIP-m-t3 | ADP<br>AMP<br>ATP<br>L=Aspartate<br>L=Glutamate<br>L=Glutamine<br>Pyrophosphate | L=2=Aminoadipate_6=semialdehyde<br>NADP+<br>NADPH | ADPM<br>AMPM<br>ATPM<br>PyrophosphateM |
|   | SC-NIP-m-t4 | ADP<br>AMP<br>ATP<br>Pyrophosphate<br>UTP | NADP+<br>NADPH | ADPM<br>AMPM<br>ATPM<br>PyrophosphateM |
|   | SC-NIP-m-t5 | ADP<br>AMP<br>ATP<br>Pyrophosphate | NADP+<br>NADPH | ADPM<br>AMPM<br>ATPM<br>PyrophosphateM |
| 3 | SC-NIP-m-t5 | ADP<br>AMP<br>ATP<br>Orthophosphate<br>Pyrophosphate | H2O2<br>NADP+<br>NADPH<br>Oxygen | ADPM<br>AMPM<br>ATPM<br>OrthophosphateM<br>PyrophosphateM |

**Table A1:** Lists of metabolites in the three heaviest components reported in Table 4.